\documentclass[a4paper, 12pt, english, reqno]{amsart}
\usepackage[ansinew]{inputenc}
\usepackage{amsmath, amssymb, latexsym, graphicx}
\usepackage[ansinew]{inputenc}
\usepackage{fontenc}
\usepackage{fancyhdr}
\usepackage{pb-diagram}
\usepackage{lastpage}
\usepackage[dvips, dvipsnames, usenames]{color}
\usepackage[stable]{footmisc}
\usepackage[backref]{hyperref}
\hypersetup{colorlinks=true, urlcolor=blue, citecolor=blue}

\newcommand{\cnoindent}{\vspace{12pt}\noindent}

\newcommand{\en}[1]{\footnote{#1.}}

\newcommand{\s}[1]{$\langle #1_{i} \rangle_{i \in \mathbb{N}}$}

\newcommand{\w}{$\omega$ }

\newcommand{\worder}{$\omega$-or\-der }
\newcommand{\wordered}{$\omega$-or\-de\-red }

\newcommand{\waorder}{$\omega^*$-or\-der }
\newcommand{\waordered}{$\omega^*$-or\-de\-red }
\newcommand{\waordering}{$\omega^*$-or\-de\-ring }
\newcommand{\waorderp}{$\omega^*$-or\-der}

\newcommand{\zas}{Z$^*$-points }
\newcommand{\za}{Z$^*$-point }
\newcommand{\zasp}{Z$^*$-points}

\begin{document}

\title {The aleph-zero or zero dichotomy}

\maketitle

\begin{center}

\vspace{-20pt}

(New and extended version with new arguments)

\vspace {20pt}

\small{

Antonio Leon Sanchez\\
I.E.S. Francisco Salinas, Salamanca, Spain\\
}
\href{http://www.interciencia.es}{http://www.interciencia.es}\\
\href{mailto:aleon@interciencia.es}{aleon@interciencia.es}

\end{center}

\begin{abstract}
This paper proves the existence of a dichotomy which being formally derived from the
topological successiveness of \waorder leads to the same absurdity of Zeno's Dichotomy
II. It also derives a contradictory result from the first Zeno's Dichotomy.
\end{abstract}

\pagestyle{myheadings}

\markboth{\small{The aleph-zero or zero dichotomy}}{\small{The aleph-zero or zero
dichotomy}}

    \pagestyle{fancy}
    \fancyhf{}
    \fancyhead[LE,RO]{\thepage}
    \fancyhead[LO]{\leftmark}
    \fancyhead[RE]{\rightmark}
    \renewcommand{\headrulewidth}{0.4pt}


\section{Introduction: Zeno's Paradoxes and Modern Science}

\noindent Zeno's Paradoxes have interested philosophers of all times\en{See
\cite{Cajori1915}, \cite{Cajori1920}, \cite{Vlastos1967}, \cite{Salmon2001},
\cite{Huggett2004}, \cite{Foster2005}, \cite{Colli2006} or \cite{Mazur2007} for
historical background} although until the middle of the XIX century they were frequently
considered as mere sophisms. From that time, and particularly through the XX century,
they became the unending source of new philosophical, mathematical and physical
discussion. Authors as Hegel James, Russell, Whitehead or Bergson\en{\cite{Hegel1995},
\cite{James1996}, \cite{Russell2001}, \cite{Whitehead1948}, \cite{Whitehead1978},
\cite{Bergson2001}, \cite{Bergson2004}} focused their attention on the challenging world
of Zeno's paradoxes. At the beginning of the second half of the XX century the pioneering
works of Black, Wisdom, Thomson, and Benacerraf\en{\cite{Black1950}, \cite{Wisdom1951},
\cite{Thomson1954}, \cite{Thomson2001}, \cite{Benacerraf1962}} introduced a new way of
discussing the possibilities of performing an actual infinity of actions in a finite time
(a performance involved in most of Zeno's paradoxes). I refer to Supertask Theory. In
fact, infinity machines, or supermachines, are our modern Achilles substitutes. A
supermachine is a theoretical device supposedly capable of performing countably many
actions in a finite interval of time. The possibilities of performing an uncountable
infinity of actions were ruled out by P. Clark and S. Read, for which they made use of
Cantor's argument on the impossibility of dividing a real interval into uncountably many
adjacent parts.\en{\cite{Clark1984}} Although supertasks have also been examined from the
perspective of nonstandard\en{\cite{McLaughlin1992}, \cite{McLaughlin1995},
\cite{Alper1997}, \cite{McLaughlin1998}} analysis, as far as I know the possibilities to
perform an hypertask along an hyperreal interval of time have not been discussed,
although finite hyperreal intervals can be divided into hypercountably many successive
infinitesimal intervals, the so called hyperfinite partitions.\en{\cite{Stroyan1997},
\cite{Goldblatt1998}, \cite{Keisler2000}, \cite{Henle2003}, etc} Supertask theory has
finally turned its attention, particularly from the last decade of the XX century,
towards the discussion of the physical plausibility of supertasks as well as on the
implications of supertasks in the physical world including relativistic and quantum
mechanics perspectives.\en{\cite{Pitowsky1990}, \cite{Laraudogoitia1996},
\cite{Laraudogoitia2001}, \cite{Salmon2001}, \cite{Grunbaum2001}, \cite{Grunbaum2001b},
\cite{Grunbaum2001c}, \cite{Laraudogoitia1997}, \cite{Laraudogoitia1998},
\cite{Earman1998}, \cite{Laraudogoitia1999}, \cite{Norton1999}, \cite{Alper1999},
\cite{Alper2000}, \cite{Laraudogoitia2002}, \cite{Weyl1949}, \cite{Hogarth1992},
\cite{Earman1993}, \cite{Earman1996}, \cite{Earman2004}, \cite{Silagadze2005}}

\cnoindent In the second half of the XX century, several solutions to some of Zeno's
paradoxes have been proposed. Most of those solutions were found in the context of new
branches of mathematics as Cantor's transfinite arithmetic, topology, measure theory and
more recently internal set theory\en{\cite{Grunbaum1955}, \cite{Grunbaum1967},
\cite{Zangari1994}, \cite{Grunbaum2001}, \cite{Grunbaum2001b}, \cite{Grunbaum2001c},
\cite{McLaughlin1992}, \cite{McLaughlin1995}} (a branch of nonstandard analysis). It is
also worth noting the solutions proposed by P. Lynds\en{\cite{Lynds2003},
\cite{Lynds2003b}} within a classical and quantum mechanics framework Some of these
solutions, however, have been contested. And in most of cases the proposed solutions do
not explain where Zeno's arguments fail. Moreover, some of the proposed solutions gave
rise to a new collection of problems so exciting as Zeno's
paradoxes.\en{\cite{Papa-Grimaldi1996}, \cite{Alper1997}, \cite{Laraudogoitia2001},
\cite{Salmon2001}, \cite{Huggett2004} \cite{Silagadze2005}}

\cnoindent The four most famous Zeno's paradoxes are usually regarded as arguments
against motion\en{\cite{Aristoteles1998}, \cite{Grunbaum1967}, \cite{Heath1981},
\cite{Cohn1994}, \cite{Salmon2001} etc} be it performed in a continuous or in a
discontinuous world. \emph{Achilles and the Tortoise} and the \emph{Dichotomy} in the
continuous case, \emph{the Stadium} and \emph{the Arrow} in the discontinuous one. The
paradoxes of the second case (together with the paradox of Plurality) are more difficult
to solve, if a solution exists after all, particularly in a quantum spacetime framework.
Most of the proposed solutions to Zeno's paradoxes are, in effect, solutions to the
paradoxes of the first group or to the second one in a dense and continuous spacetime
framework. This situation is very significant taking into account the increasing number
of contemporary physical theories suggesting the quantum nature of spacetime, as for
instance superstring theory, loop quantum gravity, quantum computation theory, or black
hole thermodynamics.\en{\cite{Green2001}, \cite{Green2004}, \cite{Veneziano2004},
\cite{FernandezBarbon2005}, \cite{Smolin2003}, \cite{Baez2003}, \cite{Smolin2004},
\cite{Susskind1997}, \cite{Bekenstein2003}, \cite{Lloyd2005}, \cite{Bekenstein2003},
\cite{Susskind1997}} Will physics (the science of change) finally meet the problem of
change (whose insolvability motivated Zeno's argument) at this quantum
level\en{\cite{Belot2001}}? Is the problem of Change really inconsistent as some
authors\en{\cite{McTaggart1908}, \cite{Mortensen2002}} have defended? These are in fact
two intriguing and still unsolved questions related to Zeno's
arguments.\en{\cite{Papa-Grimaldi1996}}

\section{Zeno's paradoxes and $\omega$-order}

\noindent Not less intriguing, though for different reasons, is the fact that one
immediately perceives when examining the contemporary discussions on Zeno's paradoxes.
Surprisingly, the Axiom of Infinity is never involved in such discussions. Zeno's
arguments have never been used to question the Axiom of Infinity, as if the existence of
actual infinite totalities were beyond any doubt.\en{\cite{Ferreiros1999}} Grünbaum, for
instance, proposed that if it were the case that from modern kinematics together with the
denseness postulate a false zenonian conclusion could be formally derived, then we would
have to replace current kinematics by other mechanical theory.\en{\cite[page
39]{Grunbaum1967}} Anything but questioning the hypothesis of the actual infinity from
which the involved topological denseness derives. And this in spite of the lack of
evidence of that hypothesis, which is even rejected by some schools of contemporary
mathematics as constructivism (among whose precursors we find scholars as Newton, Fermat
or Euler\en{\cite{Marion1998}}) and by some XX century thinkers of the intellectual
stature of Poincaré  or Wittgenstein.\en{\cite{Poincare1913}, \cite{Wittgenstein1978}}

\cnoindent In the first half of the XIX century Bernard Bolzano, and in the second one
Richard Dedekind, tried unsuccessfully\en{Their respective proofs were compatible with
the potential infinity.} to prove the existence of infinite
totalities.\en{\cite{Bolzano1993}, \cite{Dedekind1988}} For his part, G. Cantor, the
founder of transfinite mathematics, simply took it for granted the existence of such
totalities. Thus, in \S6 of his famous Beiträge (pp. 103-104 of the English translation)
we can read:

    \begin{quote}
        The first example of a transfinite set is given by the totality of finite cardinals.
    \end{quote}

\cnoindent although, as could be expected, he gave no proof on the the existence of that
totality. In accordance with his profound theological platonism, Cantor was firmly
convinced of the actual existence of complete infinite totalities.\en{\cite{Dauben1990}}
He never explicitly declared the hypothetical nature of his infinitist assertions (at
least not in his most relevant works on the transfinite.\en{\cite{Cantor1895},
\cite{Cantor1897}, \cite{Cantor1891}, \cite{Cantor1976}, \cite{Cantor1994}} He even tried
to prove the existence of actual infinities (quoted in \cite{Rucker1995}, p. 3, from
\cite{Cantor1932}, p. 404):

    \begin{quote}
        ... in truth the potential infinite has only a borrowed reality, insofar as
        potentially infinite concept always points towards a logically prior actually infinite
        concept whose existence it depends on.
    \end{quote}

\noindent Evidently this is not a formal proof but a personal belief. Cantor's infinite
totality is isomorph to the set $\mathbb{N}$ of natural numbers and then his implicit
assumption on the existence of that complete totality is equivalent to our modern Axiom
of Infinity.

\cnoindent The (assumed) infinite totality of finite cardinals led Cantor to the
essential notion of $\omega$-or\-der (Beitr\"age, p. 115 \cite{Cantor1955}):

    \begin{quote}
        By $\omega$ we understand the type of a well ordered aggregate:
        \begin{equation*}
        (e_1, e_2, \dots, e_\nu, \dots)
        \end{equation*}
        in which:
        \begin{equation*}
        e_\nu \prec e_{\nu + 1}
        \end{equation*}
        and where $\nu$ represents all finite cardinal numbers in turn.
    \end{quote}

\cnoindent Cantor then defined the notion of \emph{fundamental series of ordinals} of
which he proved the existence of a limit (Beitr\"age, Theorem \S14 I). This limit plays a
capital role in the proofs of the following 10 theorems in Beitr\"age \S 15 the last of
which is the fundamental theorem K-15 stating that every ordinal of the second class
(transfinite) is either the result of increasing by one the next smaller ordinal
(ordinals of the first kind), or the limit of a fundamental increasing sequence of
ordinals (ordinals of the second kind). Cantor construction of transfinite ordinals, from
$\omega$ to the $\epsilon$-numbers of the second number class, strongly depends on
Theorem K-15. The imposing cantorian edifice was really founded on that theorem. And that
theorem, in turn, depends on the hypothetical existence of a complete infinite totality:
that of the finite cardinal numbers (Axiom of Infinity in modern terms), which is
anything but selfevident.

\cnoindent In modern terms we say a sequence is \wordered if it has a first element and
each element has an immediate successor. Similarly, a sequence is \waordered if it has a
last element and each element has an immediate predecessor. Evidently, both type of
ordering are intimately related to Zeno's Dichotomies, although surprisingly, the
analysis of Zeno's arguments as formal consequences of \worder remains still undone. For
some unknown reasons, it seems we are not interested in analyzing the consistency of the
hypothetical existence of actual complete infinite totalities. And this in spite of the
enormous problems the actual infinity poses to experimental sciences as physics (recall
for example the problems of renormalization in elementary particle
physics\en{\cite{Feynman1988}, \cite{Hooft1991}, \cite{Green2001}, \cite{Yndurain2002},
\cite{Green2004}}). The discussion that follows is just oriented in that direction. Its
main objective is to analyze Zeno's Dichotomies I and II from the perspective of \worder
and \waorder respectively.

\section{The aleph-zero or zero dichotomy}

\noindent In what follows, and for the sake of clarity, I will consider a canonical
version of the famous Achilles' race whose logical impossibility Zeno claimed. In fact,
Achilles will be considered as a single mass point moving rightwards along the X axis,
from point -1 to point 1, at a finite velocity $v$. In the place of the uncountable and
densely ordered sequence of points within the real interval $[-1, 1]$ we will only
consider the \waordered sequence of points:
    \begin{equation}
        \dots, \ \frac{1}{2^4}, \frac{1}{2^3}, \frac{1}{2^2}, \frac{1}{2}, 1
    \end{equation}
all of which Achilles must successively traverse in order to reach point 1 from the
starting point -1. In fact, this denumerable sequence of points (\zas according to
classical Vlastos' terminology\en{\cite{Vlastos1966}}) is not densely by successively
ordered, which means that between any two successive \zas no other \za exists. In
consequence, and at a finite velocity, \zas can only be traversed in a successive way:
one after the other. Assume now Achilles is just on the point 0 at the precise instant
$t_0$. According to classic mechanics he will reach point 1 just at\en{Assuming, for
instance, that $v$ is given in kilometers per second and that the distance from 0 to 1 is
just one kilometer.} $t_1=t_0 + 1/v$. But before reaching his goal, he has to
successively traverse the controversial \zasp. We will focus our attention just on the
way Achilles performs such a traversal. For this, let $f_{z^*}(t)$ be the number of \zas
Achilles has traversed at the precise instant $t$, being $t$ any instant within the
closed real interval $[t_0, \, t_1]$. It is quite clear that $f_{z^*}(t_0)=0$ because at
$t_0$ Achilles is just on point 0. For any other instant $t$ in $[t_0, \, t_1]$ Achilles
has already passed over countably many \zasp, for if there were an instant $t$ in $[t_0,
\, t_1]$ at which Achilles were passed only over a finite number $n > 0$ of \zasp, these
$n$ \zas would have to be the impossible firsts $n$ points of an \waordered sequence of
points. So we can write:
    \begin{equation}
        f_{z^*}(t)=
        \begin{cases}
            0 &\text{if $t = t_0$}\\
            \aleph_0 &\text{if $t_0 < t \leq t_1$}
        \end{cases}
    \end{equation}
Notice $f_{z^*}(t)$ is well defined for each $t$ in $[t_0, \, t_1]$. Consequently,
$f_{z^*}$ maps the real interval $[t_0, t_1]$ into the set of two elements $\{0, \
\aleph_0\}$. In this way $f_{z^*}$ defines a clair dichotomy, the \emph{aleph-zero or
zero dichotomy},\en{Although the usual way of reading $\aleph_0$ is aleph-null -it can
also be read as aleph-zero- the original English translation by P. E. B. Jourdain of
Cantor's Beiträge was ''aleph-zero''. Section 6 is entitled ''The Smallest Transfinite
Cardinal Number Aleph-Zero. The current English edition of Cantor's Beiträge is from
1955.} regarding the numbers of \zas Achilles has traversed when moving rightward from -1
to 1 along the X axis. Accordingly, with respect to the number of the traversed \zasp,
Achilles can only exhibit two states:
    \begin{enumerate}
        \item State $A_0$: Achilles has traversed zero \zasp.
        \item State $A_{\aleph_0}$: Achilles has traversed aleph-zero \zasp.
    \end{enumerate}
\noindent Thus, Achilles directly becomes from having traversed no \za (state $A_0$) to
having traversed $\aleph_0$ of them (state $A_{\aleph_0}$). Finite intermediate states,
as $A_n$ at which he would have traversed only a finite number $n$ of \zasp, simply do no
exist. The set of states Achilles exhibits with respect to the number of traversed
Z-points is well defined and has only two elements, namely $A_0$ and $A_{\aleph_0}$. Let
us now examine the transition from $A_0$ to $A_{\aleph_0}$ under the inevitable
restriction of the above \emph{aleph-zero or zero dichotomy}. The topological
successiveness of \zas makes it impossible to traverse them other than successively. And
taking into account that between any two successive \zas a finite distance greater than 0
exists, to traverse $\aleph_0$ of those \zas -whatever they be- means to traverse a
finite distance greater than 0. This traversal, at the finite Achilles' velocity $v$, can
only be accomplished by lasting a certain amount of time necessarily greater than 0.
Achilles, therefore, has to expend a certain amount of time $\tau > 0$ in becoming
$A_{\aleph_0}$ from $A_0$. The \waordering imposes this time $\tau$ has to be
indeterminable, otherwise we would know the precise instant at which Achilles becomes
$A_{\aleph_0}$ and, consequently, we would also know the precise \za point on which he
reaches that state, and this is impossibly because in that case we would have a natural
number $n$ such that $n + 1 = \aleph_0$. The indeterminacy of $\tau$ means both the
existence of more than one alternative and the impossibility to determine the actual
alternative. Now then, indeterminable as it may be, $\tau$ must be greater than 0, and
this inevitable requirement imposed by the fact that Achilles must traverse a distance
greater than 0 at its finite velocity $v$ is incompatible with the aleph-zero or zero
dichotomy, as we will immediately see.

\cnoindent In fact, let $\tau$ be any real number greater than zero and assume the
transition from $A_0$ to $A_{\aleph_0}$ lasts a time $\tau$. Consider the real interval
$(0,\tau)$ and any instant $t \in (0, \tau)$. At $t$ Achilles state cannot be neither
$A_0$ nor $A_{\aleph_0}$. It cannot be $A_0$ because if that were the case the process of
becoming $A_{\aleph_0}$ would not have begun, and then the process of becoming
$A_{\aleph_0}$ would last a time equal or less than $\tau - t$ in the place of the
assumed $\tau$. It cannot be $A_{\aleph_0}$ because in that case the process of becoming
$A_{\aleph_0}$ would have already finished and then it would have lasted a time equal or
less than $t$ in the place of the assumed $\tau$. Now then, Achilles' state has to be
either $A_0$ or $A_{\aleph_0}$ because it is well defined along the real interval $[0,
1]$ of which $(0, \tau)$ is a proper subinterval. Consequently, and being $\tau$ any real
number, it is impossible for Achilles to become $A_{\aleph_0}$ from $A_0$ by lasting a
time greater than zero. Notice this is not a question of indeterminacy but of
impossibility: no real number greater than zero exists for the duration of Achilles'
transition from $A_0$ to $A_{\aleph_0}$. He, therefore, has to become $A_{\aleph_0}$ from
$A_0$ instantaneously. But this is impossible at his finite velocity $v$. He must,
therefore, remain $A_0$. Or in other words, he cannot begin to move. Evidently, this
conclusion is the same absurdity of Zeno's Dichotomy II, although in our case it has been
directly derived from the topological successiveness of \waorderp, which in turn derives
from assuming the existence of complete denumerable totalities (Axiom of infinity), as
Cantor proved.\en{\cite{Cantor1955}} To be complete (as the actual infinity requires) and
uncompletable (because no last -first- element completes them) could be a contradictory
attribute rather than a permissible eccentricity of both \wordered and \waordered
sequences.

\noindent The above conclusion is confirmed by the following argument.\en{A variant of
Benardete's paradox \cite{Bernadete1964}} Let us replace each \za with a mass
$Z^*$-sensor capable of emitting a visible laser beam when it is activated by any mass
passing over it. Assume the system of $Z^*$-sensors is regulated in such a way that each
sensor emits its corresponding laser beam if, and only if, it is activated and no other
laser beam is being emitted by other $Z^*$-sensor of the system. So only one laser beam
can be being emitted by the system of $Z^*$-sensors: the one corresponding to the first
activated $Z^*$-sensor, whatsoever it be. Assume Achilles performs his canonical race
from point -1 to point 1. Will any laser beam being emitted at $t_1$? Evidently not,
because it would have to be the impossible first $Z^*$-sensor of an \waordered sequence
of $Z^*$-sensors. But on the other hand, why not? Is there any reason to explain the
inevitable malfunctioning of the $Z^*$-sensors system other than the inconsistency of
assuming that it is possible to begin a sequence of discrete and successive actions
without a first action to begin? What is impossible is not motion but the actual
infiniteness of \waorderp.

\cnoindent Let us now examine Zeno's Dichotomy I under the same canonical conditions of
the above Dichotomy II. Consider again the real interval $[-1, 1]$ in the $X$ axis. Let
now \s{z} be the \wordered sequence of Z-points:
    \begin{equation}
        \frac{1}{2}, \ \frac{3}{4}, \ \frac{7}{8}, \ \dots \frac{2^{i}-1}{2^i}, \ \dots
    \end{equation}
Achilles has to traverse in his race from point -1 to point 1. Assume also we remove from
$(0, 1)$ all points except just Z-points (we would have a sort of \emph{Zeno's powder}).
In the place of a continuous race from point -1 to point 1, assume that Achilles is on
point 0 just at instant $t_0$ and then he begins to jump to $z_1$, to $z_2$, to $z_3$,
\dots, and to any point $x$ in $[1, \, 2]$ if there were no other Z-poins to jump, in
such a way that he is on each $z_i$ just at $t_i$ as a consequence of a Z-jump $j_i$,
being $t_i$ the $i$-th term of the \wordered sequence of instants \s{t} whose limit is
$t_b$.

\cnoindent The one to one correspondence $f(t_i) = z_i$ proves\en{This is the way
infinitists pretend to explain how an \wordered sequence of actions can be completed: by
pairing off two endless sequences, the one of actions the other of instants at which the
successive actions are carried out. Thus, in order to end an endless sequence of actions
we only need to pair the endless sequence of actions with the endless sequence of
instants as which they are performed, as if by pairing off two impossibilities a
possibility could result.} that at $t_b$ Achilles has completed the \wordered sequence of
Z-jumps \s{j} on the \wordered sequence of Z-points \s{z}. Thus, at $t_b$ Achilles has to
be on a point $x \geq 1$ of $[1, \, 2]$. Otherwise, if he were on a Z-point $z_i$, only a
finite number $i$ of jumps would have been performed. We now have an uncomfortable
asymmetry between the \wordered sequence of Z-jumps \s{j} and the (\w + 1)-ordered
sequence of points: the \wordered sequence of Z-points plus the last point $x$ Achilles
ends up his \wordered sequence of Z-jumps, i.e. the (\w + 1)-ordered sequence $\langle
\langle z_i \rangle_{i \in \mathbb{N}}, x\rangle$.

\cnoindent By definition Achilles is on each $z_i$ at $t_i$ as a consequence of the i-th
Z-jump $j_i$. The one to one correspondence $f(j_i)=z_i$ proves that:
    \begin{enumerate}
        \item Each z-jump $j_i$ ends on a z-point $z_i$.
        \item Consequently: No Z-jump $j_i$ makes Achilles to reach point $x$.
        \item Consequently: Achilles comes to point $x$ from no Z-point.
    \end{enumerate}
\noindent But the only actions Achilles performs from $t_0$ is the \wordered sequence of
Z-jumps \s{j} on the \wordered sequence of Z-points \s{z}. So, Achilles can only come
from a Z-point as a consequence of a Z-jump. How is then possible Achilles reaches point
$x$ at $t_b$ if none of the performed Z-jumps places him there? At this point of the
discussion, most infinitists claim that although Achilles comes to point $x$ from no
Z-point as a consequence of no Z-jump, it reaches that point at $t_b$ as a consequence of
having \emph{completed} the (uncompletable) \wordered sequence of Z-jumps \s{j}. As if
the completion of the \wordered sequence of Z-jumps \s{j} were a place one may come from.
But if the completion of the \wordered sequence of Z-jumps \s{j} means the completion of
the \wordered sequence of Z-jumps \s{s}, i.e. that each one of the countably many Z-jumps
$j_1$, $j_2$, $j_3$, \dots, \emph{and only them}, have been performed, then it is quite
clear that Achilles cannot reach point $x$ at $t_b$. Simply because no Z-jump $j_1$,
$j_2$, $j_3$, \dots ends on the point $x$. And if no Z-jump $j_1$, $j_2$, $j_3$, \dots,
end on the point $x$ and Achilles only performs Z-jumps, then he cannot end on point $x$
either. On the other hand, if the completion of the \wordered sequence of Z-jumps \s{j}
were an additional jump then we would have an (\w + 1)-ordered sequence of jumps rather
that an \wordered one. But we have proved the \wordered sequence of Z-jumps \s{j}
suffices to place Achilles on point $x$ at $t_b$. It is therefore that \wordered sequence
\s{j} which places and does not place Achilles on point $x$.

\cnoindent Achilles ends his \wordered sequence of Z-jumps on point $x$ and this final
position is unexplainable because no final jump places him there. And no final jump
places him there because no final Z-jumps exists in the \wordered sequence of Z-jumps
\s{j}. The asymmetry is quite clair: there exists a last effect (to reach point $x$) but
not a last jump causing it. Infinitist, therefore, have to make use of a mysterious last
jump by converting the completion of an \wordered sequence of jumps in a subsequent
additional jump which is different from all previously performed ones. But an \wordered
sequence of jumps plus an additional last jump is not an \wordered sequence of jumps but
an (\w + 1)-ordered one. Thus, this assumed additional jump does not solve the question,
because Achilles reaches and does not reaches point $x$ as a consequence of an \wordered
(not of an (\w + 1)-ordered) sequence of jumps.

\cnoindent As in the case of Dichotomy II, assume that each Z-point $z_i$ is provided
with a mass Z-sensor, being the system of Z-sensors regulated in such a way that, once a
Z-sensor is activated, it will be emitting its corresponding laser beam until other
Z-sensor be activated and emits its own laser beam. In consequence, once the system is
activated there will always be a Z-beam being emitted: the one corresponding to the last
activated Z-sensor. So, once activated, it is impossible to turn off the emission of
Z-beams. Assume now Achilles performs an \wordered sequence of Z-jumps on the \wordered
sequence of Z-sensorized Z-points. For the same reasons above, Achilles completes this
\wordered sequence of Z-jumps at $t_b$. And now the question is: will any laser beam
being emitted at $t_b$? According to the functioning of the Z-sensors system, once
Achilles activates the first Z-sensor by Z-jumping on the first Z-point, it is impossible
to turn off the emission of Z-beams. So, at $t_b$ a Z-beam has to be being emitted.
Although, on the other hand, no Z-beam can be being emitted at $t_b$ because, if Achilles
has completed his uncompletable \wordered sequence of Z-jumps, that Z-beam would have to
be being emitted by the impossible last Z-sensor of the \wordered sequence of Z-sensors.
Thus, if Achilles has completed the uncompletable \wordered sequence of Z-jumps, a laser
beam will and will not be being emitted by the system of Z-sensors. This seems rather
contradictory.

\providecommand{\bysame}{\leavevmode\hbox to3em{\hrulefill}\thinspace}
\providecommand{\MR}{\relax\ifhmode\unskip\space\fi MR }
\providecommand{\MRhref}[2]{%
  \href{http://www.ams.org/mathscinet-getitem?mr=#1}{#2}
} \providecommand{\href}[2]{#2}

\end{document}